\documentclass[12pt]{article}
\usepackage{amssymb,amscd}
\usepackage{graphics}

\usepackage[dvipdfm]{graphicx}
\newtheorem{conj}{Conjecture}[section]
\newtheorem{theo}[conj]{Theorem}

\newtheorem{prop}[conj]{Proposition}

\begin{document}
\title{\Large Calculating Euler-Poincar\' e Characteristic Inductively}

\author{\Large Milosav M. Marjanovi\'c\\
{\small Serbian Academy of Sciences and Arts}\\[-2mm]
{\small Kneza Mihaila 35, 11001 Belgrade, Serbia}
\\[-2mm] {\small milomar$@$beotel.rs}}
\maketitle \vskip 2cm

\begin{abstract}
Let $((X_i,A_i)), i=0,1,...,n$ be a sequence of pairs of
topological spaces and $(Y_j), j=1,...,n$ a sequence of
topological spaces. We suppose that all spaces $X_i, i=0,1,...,n$
and $Y_j, j=1,...,n$, taken together, are mutually disjoint and
let $Y$ be the disjoint topological sum of the spaces $X_i,
i=0,1,...,n$ and $Y_j\times [j-1,j], j=1,...,n$. Let the mappings
$$f_0 : Y_1\times \{0\}\rightarrow A_0, \underline {f_i} : Y_i\times
\{i\}\rightarrow A_i, \overline {f_i} : Y_{i+1}\times
\{i\}\rightarrow A_i, f_n : Y_n\times \{n\}\rightarrow A_n$$

\noindent be continuous and onto. When for each $i$ and each $a\in
A_i$, the point $a$ is identified with all points of the sets
$\underline {f_i}^{-1}(a)$ and $\overline {f_i}^{-1}(a)$,
($\underline {f_0}=\overline {f_0}=f_0, \underline {f_n}=\overline
{f_n}=f_n$) and all other points of $Y$ with themselves, then a
quotient space $X$ is obtained. The spaces $X_i$ and $Y_i\times
\{t\}, t\in (i-1,i)$ are homeomorphic to their embedded copies
$\overline{X_i}$ and $\overline{Y_i}(t)$, called fibers of $X$ and
these fibers make an ordered decomposition of $X$, called its
fibrous decomposition and denoted by $X_0(Y_1)...(Y_n)X_n$.

We call a finite space $0$-fibrous and, proceeding inductively, we
call a space $X$ $m$-fibrous when $X$ has a fibrous decomposition
each fiber of which is $k$-fibrous for some $k$ less than $m$. We
prove that for an $m$-fibrous space $X$ its Euler-Poincar\' e
characteristic is defined and if $X_0(Y_1)...(Y_n)X_n$ is its
fibrous decomposition, then
$$\chi (X)=\chi (X_0)-\chi (Y_1)+\cdots -\chi (Y_n)+\chi (X_n).$$

Examples of calculation of E-P characteristic of a number of
spaces is given without any use of their combinatorial structures.
But when $K$ is a finite $n$-dimensional CW complex, then we find
that $K$ is an $n$-fibrous space whose E-P characteristic is $\sum
(-1)^i\alpha_i$, where $\alpha_i$ is the number of $i$-cells of
$K$.
\end{abstract}

\section{Introduction}

Teaching a course of didactics of mathematics for preservice
primary school teachers, I included a number of lectures on
topological, projective and metric properties, experienced when
visible shapes are observed. According to J. Piaget a preschool
child forms some spontaneous intuitive concepts related to the
shape of things in his/her surroundings, following the order:
topological-projective-metric. To make these ideas based on a
solid ground, I employed some basic mathematics that these
students know from their secondary school (and the mathematics
course usually scheduled for the teacher training faculties).

Besides some intuitively easy to describe topological properties,
I also included calculation of Euler-Poincar\' e (abbreviated E-P)
characteristic decomposing lines into running sets of points and
surfaces into running sets of lines (see examples 1, 2, and 4, in
the section 3. Examples of this paper and those in the paper
\cite{m}). My students were particularly excited to see a shape be
heavily distorted and still preserving its E-P characteristic.

In search of some sources where calculation of E-P characteristic
would be treated inductively, I came across some interesting
papers (for instance, \cite{cgr}, \cite{fls}, \cite{v}), but none
corroborating my unfounded method from \cite{m}. Thus, I write
this short note for my sins.

At the end, we add that the objectives of this note are more
modest than those of the papers \cite{cgr}, \cite{fls} and
\cite{v} and that we have evidently been motivated by a basic
insight provided by Morse theory.

\section{Fibrous decompositions of spaces}

All spaces that we consider are supposed to be Hausdorff. Given a
pair of spaces $(X, A)$ and a space $Y$, we suppose that the
spaces $X$ and $Y\times I$, ($I = [0,1]$) are disjoint and that $W
= X\oplus (Y\times I)$ is their disjoint topological sum. Let the
mapping $f : Y\times \{0\}\rightarrow A$ be continuous and onto.

Identifying each point $x\in A$ with all points of $f^{-1}(x)$, a
quotient space is obtained which will be denoted by
$\overline{W}$. We will also say that $\overline{W}$ is obtained
from $W$ joining together $X$ and $Y\times I$ by the mapping $f$.
The mapping $x\mapsto [x]$, which maps each point $x$ of $X$ onto
its equivalence class in $\overline{W}$ is a homeomorphism and the
homeomorphic copies of $X$ and $A$ in $\overline{W}$ will be
denoted by $\overline{X}$ and $\overline{A}$ respectively. Now we
prove a statement which will be used later in some proofs that
follow.

\begin{prop}
\label{prop1} Let $\overline{W}$ be the space obtained by joining
together $X$ and $Y\times I$ by the mapping $f$. Then, the space
$\overline{X}$ is a strong deformation retract of the space
$\overline{W}$.
\end{prop}

\medskip\noindent
{\bf Proof.} Let $\alpha : W\times I\rightarrow W$ be given by
$\alpha (x,u) = x$ for each $x\in X$ and each $u\in I$ and let
$\alpha ((y,t),u) = (y,t(1-u))$, for each $y\in Y$ and $t\in I$
and $u\in I$. Let $p : W\rightarrow \overline{W}$ be the natural
projection and $p\times i : W\times I\rightarrow \overline{W}
\times I$ be given by $(p\times i) (w,u) = ([w],u)$. Let $H :
\overline{W}\times I\rightarrow \overline{W}$ be given by
$H([x],u) = [x]$, for each $x\in X$ and $u\in I$ and $H([(y,t)],u)
= [(y,t(1-u))]$, for each $y\in Y$, $t\in I$ and $u\in I$. Since
$p\circ \alpha = H\circ (p\times i)$ and being $p\circ \alpha$
continuous and $p\times i$ quotient, it follows that $H$ is
continuous. Thus $H$ is a strong deformation retraction and
$\overline{X}$ a strong deformation retract of $\overline{W}$.
\bigskip

Now we describe a quotient model which will be the basis for the
calculation of the E-P characteristic. Let $((X_i,A_i)),
i=0,1,...,n$ be a sequence of pairs of topological spaces and
$(Y_j), j=1,...,n$ a sequence of topological spaces. We suppose
that all spaces $X_i, i=0,1,...,n$ and $Y_j, j=1,...,n$, taken
together, are mutually disjoint and let $Y$ be the disjoint
topological sum of the spaces $X_i, i=0,1,...,n$ and $Y_j\times
[j-1,j], j=1,...,n$.

Let the mappings
$$f_0 : Y_1\times \{0\}\rightarrow A_0, \underline {f_i} : Y_i\times
\{i\}\rightarrow A_i, \overline {f_i} : Y_{i+1}\times
\{i\}\rightarrow A_i, f_n : Y_n\times \{n\}\rightarrow A_n$$

\noindent be continuous and onto, for each $i$ ($i = 1,...,n-1$).

\begin{figure}[hbt]
\centering
\includegraphics[scale=.50]{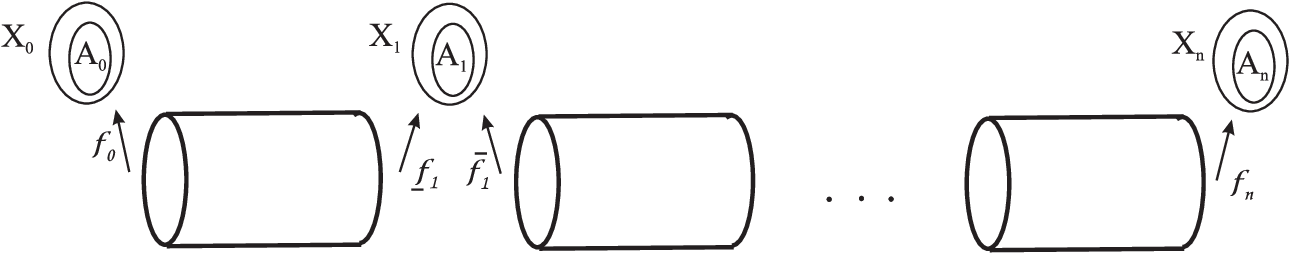}
\caption{}
\end{figure}

Let us write formally $\underline{f_0} = \overline{f_0} = f_0$ and
$\underline{f_n} = \overline{f_n} = f_n$. Now we suppose that for
each $i=0,...,n$ and each $a\in A_i$, the point $a$ is identified
with all points of the sets $\underline{f_i}^{-1}(a)$ and
$\overline{f_i}^{-1}(a)$ and all other points of $Y$ with
themselves. The quotient space which is obtained by this
identification will be denoted by $X$.

The spaces $X_i$ and $Y_i\times \{t\}, t\in (i-1,i)$ are
homeomorphic to their embedded copies $\overline{X_i}$ and
$\overline{Y_i}(t)$, called {\it fibers} of $X$ and these fibers
make an ordered decomposition of $X$, called its {\it fibrous
decomposition} and denoted by $X_0(Y_1)...(Y_n)X_n$. The number
$n$ will be called the {\it length} of the corresponding fibrous
decomposition. (When $n=0$, the decomposition reduces to $X_0$.)

Mapping each fiber $\overline{X_i}$ onto $i$ and each
$\overline{Y_i}(t)$, $t\in (i-1,i)$ onto $t$, a function $\varphi
: X\rightarrow [0,n]$ is defined and $\varphi^{-1}([0, k]), (k <
n)$ is a subspace of $X$ whose fibrous decomposition is determined
by the subsequences $(X_0,A_0),...,(X_k,A_k), Y_1,...,Y_k$ and by
the subset $\underline{f_i}, \overline{f_i}, i=0,...,k$ of the
corresponding set of mappings. We call $\varphi$ a {\it function
associated with the given fibrous decomposition}. (As it is seen,
we use a terminology to avoid confusion with the existing one
based on the morpheme "fiber").

A finite space $X$ will be called $0$-{\it fibrous} and the space
$X$ itself will be considered as its own fibrous decomposition.
Proceeding inductively, we call a topological space $m$-{\it
fibrous} when it has a fibrous decomposition each fiber of which
is $k$-fibrous for some $k\leq m-1$. Now we are ready to prove the
statement that follows.

\begin{theo}
\label{theo} Let $X$ be an $m$-fibrous space having its fibrous
decomposition given by the sequences $(X_0,A_0),(X_1,A_1),...,
(X_n,A_n)$ and $Y_1,...,Y_n$ and the set of mappings
$\underline{f_i}, \overline{f_i}, i = 0,1,...,n$. Then, the
Euler-Poincar\' e characteristic is defined for all spaces $X_0,
X_1,...,X_n$ and $Y_1,...,Y_n$ as well as for $X$ and

$$\chi (X) = \chi (X_0)-\chi (Y_1)+\chi (X_1)-\cdots +
\chi (X_{n-1})-\chi (Y_n)+\chi (X_n).$$
\end{theo}

\medskip\noindent
{\bf Proof.} The statement is trivially true when $m = 0$. Let us
suppose that it is true for all spaces which are $k$-fibrous for
some $k\leq m - 1$. Let $X$ be an $m$-fibrous space. When $X$ has
a fibrous decomposition of the length $0$, then $X$ is $k$-fibrous
for some $k\leq m-1$ and the statement is true. Let us suppose it
is true for all $m$-fibrous spaces having a fibrous decomposition
of the length less than $n$. Let $X$ be an $m$-fibrous space and
let $(X_0,A_0),(X_1,A_1),...,(X_n,A_n)$ and $Y_1,...,Y_n$ be
sequences which, together with the set of mappings
$\underline{f_i}, \overline{f_i}, i=0,1,...,n$ determine a fibrous
decomposition of $X$. According to the definition of $m$-fibrous
spaces, the spaces $X_0,X_1,...,X_n$ and $Y_0,...,Y_n$ are
$k$-fibrous for some $k\leq m-1$ and according to the inductive
hypothesis on $m$, the E-P characteristic is defined for them.

Let $\varphi$ be the function associated with the fibrous
decomposition of $X$. Modifying slightly Proposition \ref{prop1},
it is easily proved that $\varphi^{-1}([0,n-1])$ is a strong
deformation retract of $\varphi^{-1}([0,n-1/3))$ as it is
$\varphi^{-1}(n)$ of $\varphi^{-1}((n-2/3,n])$. As we have already
noticed it, the E-P characteristic is defined for $\varphi^{-1}(n)
\approx X_n$ and, according to the induction hypothesis on $n$, it
is also defined for $\varphi^{-1}([0,n-1])$. From the following
homotopy equivalences

$$\varphi^{-1}([0,n-1])\simeq \varphi^{-1}([0,n-1/3)),
\varphi^{-1}(n)\simeq \varphi^{-1}((n-2/3,n])$$

\noindent it follows that E-P characteristic is also defined for
the spaces $\varphi^{-1}([0,n-1/3))$ and
$\varphi^{-1}((n-2/3,n])$. Being these two spaces open in $X$, we
see that the triad

$$(X,\varphi^{-1}([0,n-1/3)),\varphi^{-1}((n-2/3,n]))$$

\noindent satisfies the excision property. Using now a very well
known property of E-P characteristic (see, for example, \cite{d}),
we can write

$$\chi (X) =\chi (\varphi^{-1}([0,n-1/3)))+\chi (\varphi^{-1}
(n-2/3,n]))-\chi (\varphi^{-1}((n-2/3,n-1/3)))$$.

Since $\varphi^{-1}((n-2/3,n-1/3))\simeq Y_n$ and using previously
established homotopy equivalences, we have

$$\chi (X) =\chi (\varphi^{-1}([0,n-1]))+\chi (\varphi^{-1}(n))-\chi
(Y_n).$$

Using the induction hypothesis on $n$, we replace $\chi
(\varphi^{-1}([0,n-1]))$ by the corresponding alternating sum,
obtaining so the following equality

$$\chi (X)=(\chi (X_0)-\chi (Y_1)+\cdots +\chi (X_{n-1}))+\chi (X_n)-\chi
(Y_n).$$

Thus, we have proved all conclusions of Theorem \ref{theo}.

\section{Examples}

First we prove a simple statement which is often useful when a
quotient model of a space is replaced with a more convenient one
defined on its subspace.

Let $X$ be a topological space and $\rho$ an equivalence relation
on $X$. Let $A$ be a subset of $X$ such that for each $x\in X$,
$[x]\cap A\neq \emptyset $, where $[x]$ is the equivalence class
of $x$. Taking $A$ with its relative topology, the induced
equivalence relation $\rho_A$ on $A$ determines the quotient space
$A/\rho_A$. The model of the quotient space $A/\rho_A$ is simpler
than that of $X/\rho$ and it is of some interest to know under
which conditions these two quotient spaces are homeomorphic.

Let $p : X\rightarrow X/\rho$ and $p_A : A\rightarrow A/\rho_A$ be
the natural projections and $i : A\rightarrow X$ and $j : A/\rho_A
\rightarrow X/\rho$ the inclusions. Then, $p\circ i = j\circ p_A$.
Being $p\circ i$ continuous and $p_A$ quotient, it follows that
$j$ is also continuous. Being $j$ 1-1 and onto, $j^{-1}$ is also
defined and we are looking for the conditions under which it is
also continuous (and therefore, when $X\setminus A$ can be cut out
from $X$).

\begin{prop} Let $(X,\rho )$ be a topological space with an
equivalence relation on $X$. Let $A\subset X$ be such that for
each $x\in X$, $[x]\cap A\neq \emptyset $ and let $\rho_A$ be the
induced relation on $A$. If one of the following conditions

      (i) $A$ is open and $p : X\rightarrow X/\rho$ is open,

     (ii)  $A$ is closed and $p : X\rightarrow X/\rho$ is closed,

    (iii)  $A$ is compact and $X/\rho$ is Hausdorff

\noindent holds true, then $A/\rho_A \approx X/\rho$.
\end{prop}

\medskip\noindent
{\bf Proof.} Under (i), ((ii)) the mapping $p\circ i$ is open
(closed) and onto. Hence, $p\circ i$ is quotient. From
$j^{-1}\circ (p\circ i) = p_A$, it follows that $j^{-1}$ is
continuous.

Under (iii), $p\circ i$ maps closed (compact) subsets of $A$ onto
closed (compact) subsets of $X/\rho$. Hence, $p\circ i$ is
quotient, what implies that $j^{-1}$ is continuous.
\bigskip

When we use $X_0(Y_1)X_1...(Y_n)X_n$ to denote a fibrous
decomposition, then the fibers $\overline{X_i}$ are called {\it
transitional} and $\overline{Y_i}(t)$ {\it running}.
\medskip

{\bf Example 1.} Let $X$ is a finite space having $n$ points.
Then, $\chi (X)= n$.
\medskip

{\bf Example 2.} Let $X$ is a rosette of $n$ circles. Then, $\chi
(X) = 1-n$.
\smallskip

For a single copy of $S^1$, $\chi (S^1) = 1-2+1 = 0$. Let us
suppose that E-P characteristic of a rosette of $n-1$ circles is
$1-(n-1)$. Let $X_0$ be the subspace of $X$ that consists of $n-1$
circles.

\begin{figure}[hbt]
\centering
\includegraphics[scale=.50]{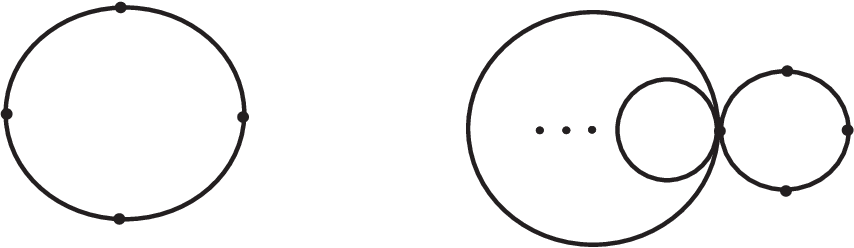}
\caption{}
\end{figure}

A fibrous decomposition of $X$ is $X_0(2p)p$, where $p$ and $2p$
denote a one point and a two point space, respectively. Thus, we
have

$$\chi (X) = \chi (X_0)-2+1 = 1-(n-1)-2+1 = 1-n.$$

Following a similar proof, it is easy to see that for the space
$X$ which is the sequence of touching circles $(x - 2k)^2+y^2 = 1,
k = 0,1,...,n-1$, $\chi (X) = 1-n$.
\medskip

{\bf Example 3.} When the boundary of an $n$-ball collapses to a
point, an $n$-sphere $S^n$ is obtained and one of its fibrous
decompositions is $p(S^{n-1})p$. From $\chi (S^n) = 1-\chi
(S^{n-1})+1$, starting with $\chi (S^0) = 2$ and applying
induction, one finds that $\chi (S^n)$ is $2$ for $n$ even and $0$
for $n$ odd.
\medskip

{\bf Example 4.} The case of surfaces (both orientable and non-orientable) deserves a special attention. Here is a brief description
of the corresponding fibrous decompositions  while the detailed presentation of examples is postponed for Section \ref{sec:surfaces}.

Let $M_g$ be the surface which is $2$-sphere with
$g$ holes (Fig. 3).

\begin{figure}[hbt]
\centering
\includegraphics[scale=.50]{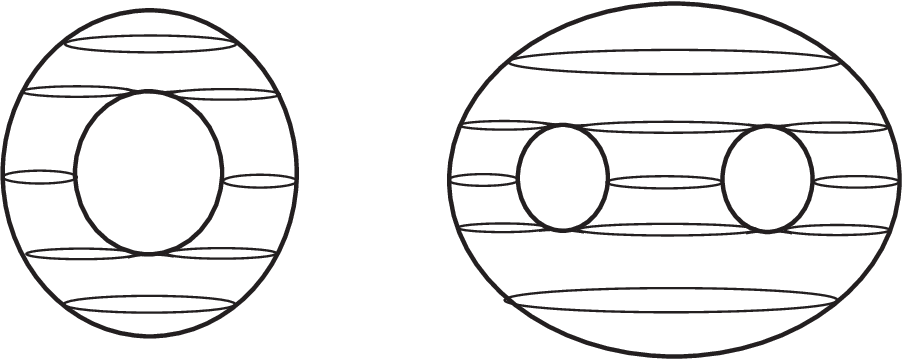}
\caption{}
\end{figure}

A fibrous decomposition of $M_g$ is $p(S^1)\dot
S^1_{g+1}((g+1)S^1)\dot S^1_{g+1}(S^1)p$, where $\dot S^1_{g+1}$
is the sequence of $g+1$ touching circles (Example 2) and
$(g+1)S^1$ disjoint topological sum of $g+1$ circles. Hence,

$$\chi (M_g) = 1-0+(-g)-0+(-g)-0+1 = 2-2g.$$

When $M_g$ is given as the quotient space obtained by
identification of arcs of the boundary of a disc, following the
command $\alpha_1\beta_1\alpha_1^{-1}\beta_1^{-1}...
\alpha_g\beta_g\alpha_g^{-1}\beta_g^{-1}$, then $M_g$ has a
fibrous decomposition $p(S^1)\tilde{S^1}_{2g}$, where
$\tilde{S^1}_{2g}$ is a rosette of $2g$ circles. Now, we have

$$\chi (M_g) = 1-0+(1-2g) = 2-2g.$$

Identifying the opposite points of boundary circles of a disk, a
disk with a circular hole, a disk with $2$ circular holes, etc.
the surfaces $N_1$ (projective plane), $N_2$ (Klein bottle),
$N_3$, etc. are obtained.

\begin{figure}[hbt]
\centering
\includegraphics[scale=.50]{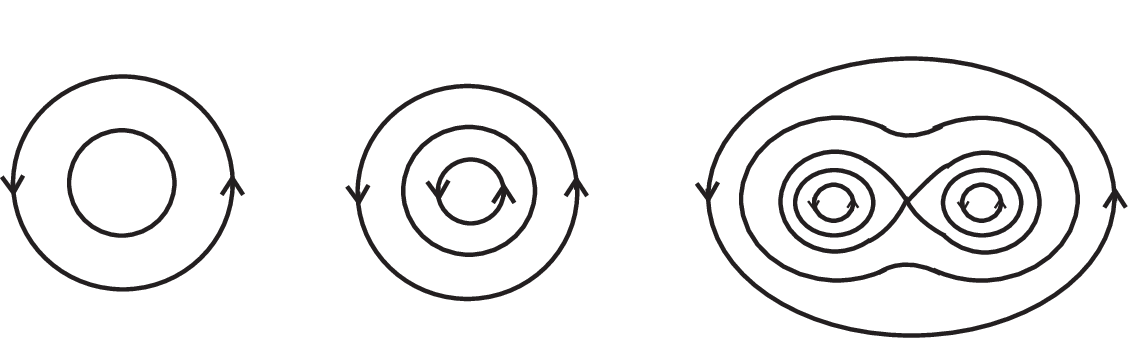}
\caption{}
\end{figure}

Their fibrous decompositions are: $S^1(S^1)p$, $S^1(S^1)S^1$,
$S^1(S^1)\dot S^1_2(2S^1)2S^1$, (where $\dot S^1_2$ are two
touching circles), etc. Their E-P characteristics are
$$\chi (N_1)=0-0+1=1, \chi (N_2)=0-0+0=0, \chi
(N_3)=0-0+(-1)-0+0=-1,$$

\noindent etc. In the case of the surface $N_h$, one of its
fibrous decompositions is

$$S^1(S^1)\dot S^1_{h-1}((h-1)S^1)(h-1)S^1,$$

\noindent where $\dot S^1_{h-1}$ is a sequence of $h-1$ touching
circles. Hence,

$$\chi (N_h) = 0-0+(1-(h-1))-0+0 = 2-h.$$

Taking $N_h$ as a quotient space obtained by identification of
boundary points of a disk, following the command
$\alpha_1^2...\alpha_h^2$, one of fibrous decompositions of $N_h$
is $p(S^1)\tilde{S^1}_h$, where $\tilde{S^1}_h$ is the rosette of
$h$ circles. Hence, $\chi (N_h) = 1-0+(1-h) = 2-h$.
\medskip

{\bf Example 5.} Let $\mathbb{R}P^n$ be $n$-dimensional real
projective space obtained when antipodal points of the boundary of
an $n$-ball $B^n$ are identified. The quotient space obtained by
this identification on the boundary $S^{n-1}$ of $B^n$ is
$(n-1)$-dimensional real projective space, what is easily seen
when open south hemisphere is cut out and Proposition 3 applied.
Thus, a fibrous decomposition of $\mathbb{R}P^n$ is
$p(S^{n-1})\mathbb{R}P^{n-1}$ and

$$\chi (\mathbb{R}P^n) = 1-\chi (S^{n-1})+\chi
(\mathbb{R}P^{n-1}).$$

Starting with $\chi (S^0) = 2$ and $\chi (\mathbb{R}P^0) = 1$, it
follows that $\chi (\mathbb{R}P^n)$ is $1$ for $n$ even and $0$
for $n$ odd.
\medskip

{\bf Example 6.} According to the way how $n$-dimensional dunce
hat $D^n$ is obtained from $n$-dimensional simplex by the
identification of points on its boundary (see \cite{ams}), the
space $D^n$ has a fibrous decomposition $p(S^{n-1})D^{n-1}$,
whence

$$\chi (D^n) = 1-\chi (S^{n-1})+\chi (D^{n-1}).$$

From this equality it easily follows that $\chi (D^n)$ is $1$ for
$n$ even and $0$ for $n$ odd.
\medskip

{\bf Example 7.} Identifying the opposite $2$-faces of the cube
$I^3 = [0,1]\times [0,1]\times [0,1]$, pairs of points $(0,y,z),
(1,y,z); (x,0,z), (x,1,z) and (x,y,0), (x,y,1)$ are identified and
a $3$-dimensional torus $T^3$ is obtained. Now we calculate
directly E-P characteristic of $T^3$ (and as it is a very well
known fact, each manifold of odd dimension has its E-P
characteristic equal $0$).

\begin{figure}[hbt]
\centering
\includegraphics[scale=.50]{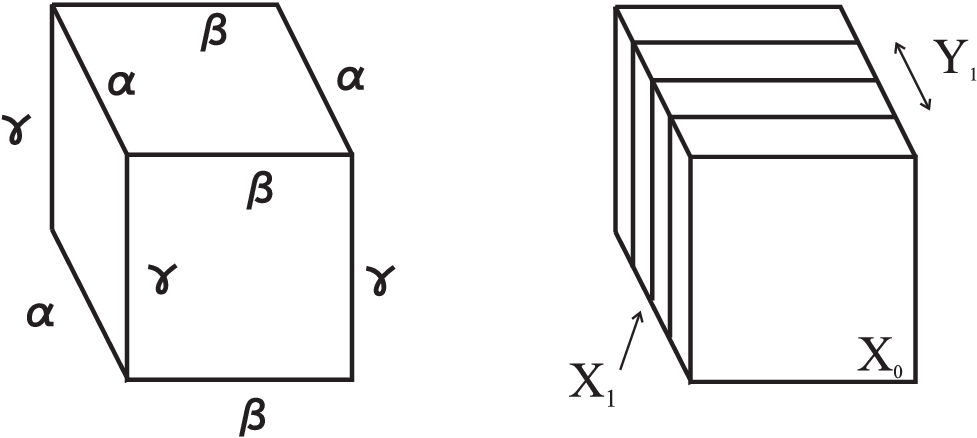}
\caption{}
\label{fig:slika-5}
\end{figure}

An obvious fibrous decomposition of $T^3$ is $T^2(2T^2)T^2$ (see
Fig.\ \ref{fig:slika-5}) and thus, $\chi (T^3) = 0-0+0 = 0$.

{\bf Example 8.} Let $K$ be a finite $n$-dimensional CW-complex.
Starting with the centers of $n$-cells, a fibrous decomposition of
$K$ is $\alpha_np(\alpha_nS^{n-1})K^{n-1}$, where $\alpha_n$ is
the number of $n$-cells of $K$ and $K^{n-1}$ is $(n-1)$-skeleton
of $K$. From that decomposition, applying induction, it easily
follows that $K$ is an $n$-fibrous space and from $\chi (K) =
\alpha_n -\alpha_n\chi (S^{n-1})+\chi (K^{n-1})$ that $\chi (K) =
\sum (-1)^i\alpha_i$, where $\alpha_i$ is the number of $i$-cells
of $K$.

\section{Euler characteristic of surfaces}\label{sec:surfaces}

More technical aspects of fibrous decompositions of surfaces are presented in Example 4, in the previous section.  Here we emphasize that
the  models of $2$-surfaces, which represent their decompositions into lines, serve for easy calculation of
the Euler-Poincar\'{e} characteristic and at the same time they enrich our geometric imagination.

\medskip
Recall that traditionally the Euler characteristic is defined  for topological spa\-ces which are homeomorphic to (finite) simplicial complexes and calculated as $\chi(X) = \sum f_i$, where $f_i$ is the number of $i$-dimensional faces of the  complex.

\medskip

Here the calculation proceeds without triangulation and, speaking figuratively, just by decomposing curves into points, surfaces into curves, bodies into surfaces, etc.

\begin{figure}
    \begin{center}
    \includegraphics[scale=.55]{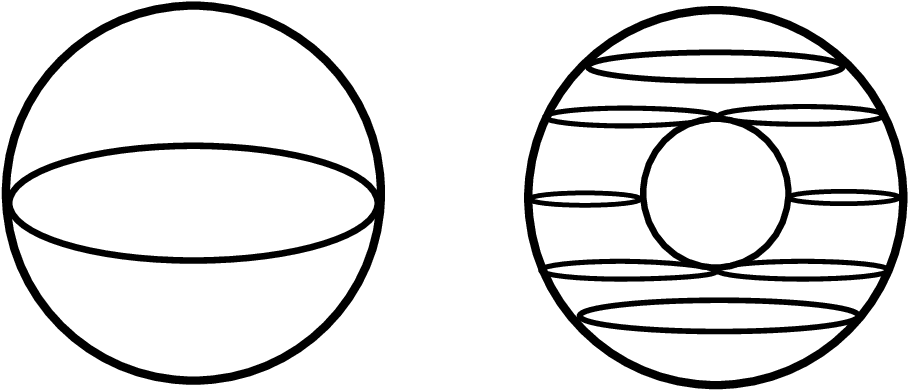}\\
    \caption{The sphere $M_0$ and the torus $M_1$}\label{fig:M0-M1}
    \end{center}
 \end{figure}

\begin{figure}
    \begin{center}
    \includegraphics[scale=.55]{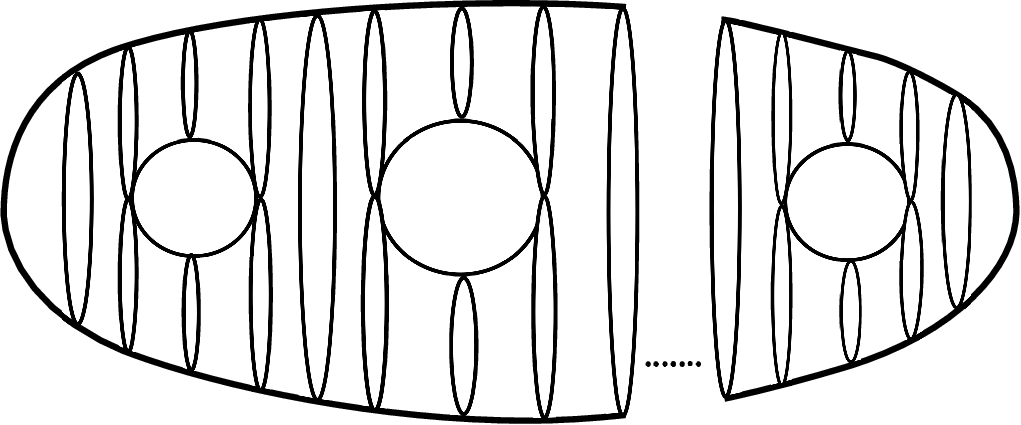}\\
    \caption{The $2$-sphere with $g$-holes}\label{fig:Mg-2}
    \end{center}
\end{figure}

\medskip\noindent
(a)  Orientable surfaces $M_g$ (spheres  with $g$ holes). The $2$-sphere and the torus
($\chi(M_0) =2, \chi(M_1) =0 $ ) are exhibited in Figure~\ref{fig:M0-M1}. The general case ($\chi(M_g) = 2-2g$) is depicted in Figure~\ref{fig:Mg-2}.

\begin{figure}
    \begin{center}
    \includegraphics[scale=.55]{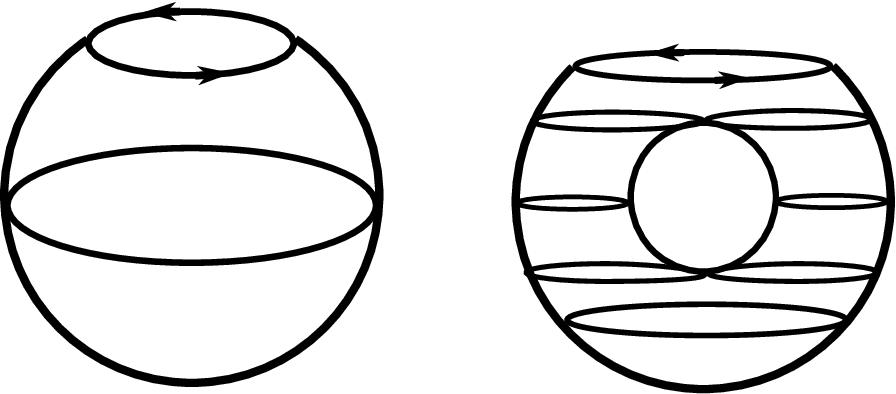}\\
    \caption{Projective plane $M_0'$ and the non-orientable surface $M_1'$}\label{fig:M0-M1-2}
    \end{center}
\end{figure}

\begin{figure}
    \begin{center}
    \includegraphics[scale=.55]{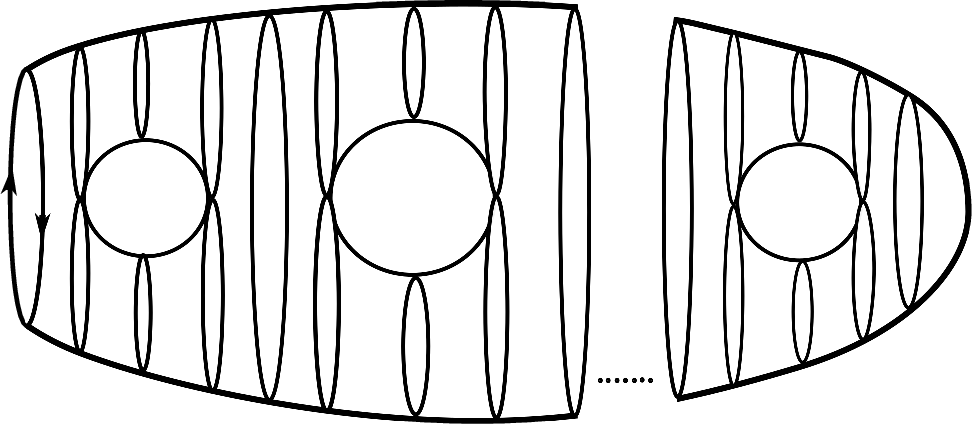}\\
    \caption{Non-orientable surface $M_g'$ with $g$-holes}\label{fig:Mg-3}
    \end{center}
\end{figure}

\begin{figure}
    \begin{center}
    \includegraphics[scale=.5]{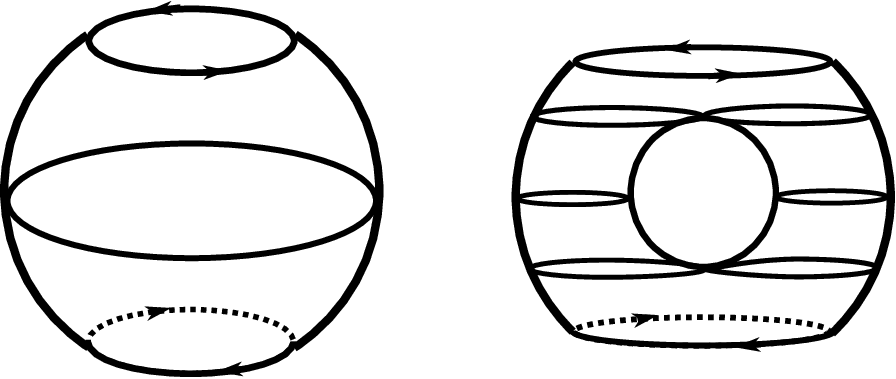}\\
    \caption{Klein bottle $M_0''$ and the non-orientable surface $M_1''$}\label{fig:M0-M1-4}
    \end{center}
\end{figure}

\begin{figure}[htb]
    \begin{center}
    \includegraphics[scale=.5]{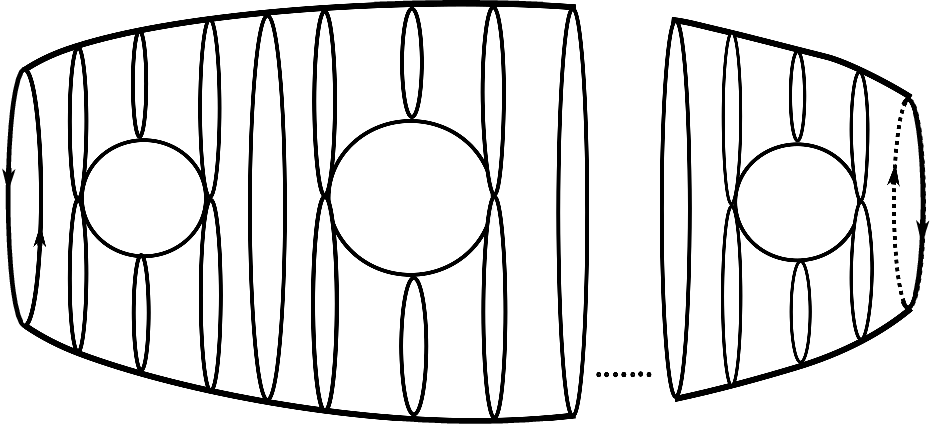}\\
    \caption{Non-orientable surface $M_g''$ with $g$-holes}\label{fig:Mg-4}
    \end{center}
\end{figure}
\noindent
(b)  Non-orientable surfaces $M_g'$ arise when a spherical cap is cut off along a circle whose diametrally
opposite points are identified. The projective plane $M_0'$ and the non-orientable surface $M_1'$
($\chi(M_0') =1, \chi(M_1') =-1 $ ) are exhibited in Figure~\ref{fig:M0-M1-2}. The general case ($\chi(M_g') = 1-2g$) is
depicted in Figure~\ref{fig:Mg-3}.

\medskip\noindent
(c) Non-orientable surfaces $M_g''$ arise when two opposite spherical caps are removed and the corresponding boundary circles are identified. The Klein bottle $M_0''$ and the non-orientable surface $M_1''$ ($\chi(M_0'') =0, \chi(M_1'') =-2$) are exhibited in Figure~\ref{fig:M0-M1-4}. The general case ($\chi(M_g'') = -2g$) is
depicted in Figure~\ref{fig:Mg-4}.

\end{document}